\documentclass[11pt]{article}
\usepackage[dvips]{graphics}
\input{amssym.def}
\input{amssym.tex}

\title{\Large\bf Riemann-Hilbert problem on an elliptic surface and a uniformly stressed inclusion embedded into a half-plane
subjected to antiplane strain}
\author{\bf Y.A.\ Antipov\\ 
Department of Mathematics, Louisiana State University\\
Baton Rouge LA 70803, U.S.A.}

\date{}

\newcommand{\I}{\mathop{\rm Im}\nolimits}
\newcommand{\R}{\mathop{\rm Re}\nolimits}
\newcommand{\const}{\mbox{const}}

\newcommand{\Md}{\partial}

\newcommand{\ov}[1]{\overline{#1}}

\newcommand{\Gg}{\gamma}
\newcommand{\Gc}{\chi}

\newcommand{\Gk}{\kappa}

\newcommand{\Gl}{\lambda}

\newcommand{\Gs}{\sigma}

\newcommand{\Gj}{\tau}

\newcommand{\Go}{\omega}
\newcommand{\Gx}{\xi}

\newcommand{\Gz}{\zeta}

\newcommand{\GF}{\Phi}

\newcommand{\GO}{\Omega}

\newcommand{\CD}{{\cal D}}

\newcommand{\CF}{{\cal F}}

\newcommand{\CI}{{\cal I}}

\newcommand{\CL}{{\cal L}}

\newcommand{\CR}{{\cal R}}

\newcommand{\beq}{\begin{equation}}
\newcommand{\eeq}{\end{equation}}

\newcommand{\fr}{\frac}

\begin{document}
\maketitle

\begin{abstract}

An inverse problem of elasticity of $n$ elastic inclusions embedded into an elastic half-plane is analyzed.
The boundary of the half-plane is free of traction. The half-plane and the inclusions are subjected to antiplane shear,
and the conditions of ideal contact hold in the interfaces between the inclusions and the half-plane.
The shapes of the inclusions  are not prescribed and have to be determined by enforcing uniform stresses
inside the inclusions. The method of conformal mappings from a slit domain onto the 
$(n+1)$-connected physical domain is worked out.  It is shown that to recover the map
and therefore the inclusions shapes, one needs to solve a vector Riemann-Hilbert problem on a 
genus-$n$ hyperelliptic surface. In a particular case of loading of a single  inclusion  in a half-plane, the problem is equivalent to two scalar Riemann-Hilbert problems
on two slits on an elliptic surface. 
In addition to three parameters of the model the conformal map possesses a free geometric parameter.
Results of numerical tests which show the impact of these parameters on the  inclusion shape are presented.

 \end{abstract}

\setcounter{equation}{0}

\section{Introduction}

Methods of conformal mappings have numerous applications in model problems of continuum mechanics.
This particularly concerns inverse problems for multiply connected domains arising in fluid mechanics
and elasticity. The former includes free boundary problems of Hele-Show and Muscat flow \cite{gus}, \cite{cro},
supercavitating flow \cite{antsil}, vortex dynamics \cite{chr}, \cite{antzem}.
The study of inverse elastic problems of determination of cavities and inclusions profiles  was initiated in \cite{esh}.  One of the first such problems, the Cherepanov
model  \cite{che},
concerns an elastic plane with $n$ cavities subjected to constant normal and tangential traction components on the their boundary 
when
the holes shapes are to be determined by enforcing constant tangential normal stresses on the boundary.
In the symmetric case of two cavities their shapes were recovered \cite{che} by applying a conformal mapping
and solving two Schwarz problems on two silts. 
Methods of integral equations \cite{vig}, the Riemann-Hilbert problems
on a hyperelliptic surface of genus $n-1$ \cite{ant1}, and the Schottky-Klein prime function \cite{mar} were
developed to generalize the solution \cite{che} to the case of any number of holes and analyze the properties 
of the conformal mappings employed and the solution derived. The method of the Riemann-Hilbert problem
on a hyperelliptic surface was worked out \cite{obn} to construct a meromorphic solution to the elastic-plastic antiplane
model for a multiply-connected domain \cite{che2}.

Inverse elastic problems of antiplane strain in multiply connected domains have been attracting
attention of many researches. The case of two symmetric inclusions was analyzed in \cite{kan}  by using the Weierstrass
zeta function and the Schwarz-Christoffel formula.  Another approach to the problem
of two finite inclusions not necessarily symmetric was proposed in \cite{wan1}. The method was designed for a specific case of uniform
stress distribution inside the inclusions
and
 based on the  Laurent
series representation of a conformal map from an annulus to the exterior of two  inclusions. 
This approach was recently applied  \cite{wan2} to the case when one of the uniformly stressed inclusions was finite, while the second
one was a semi-infinite body.
 Numerical solutions for multiply-connected
domains were obtained  by the method of
finite elements in \cite{liu} and by Faber series with the coefficients determined by nonlinear systems in \cite{dai}.

Two methods of conformal mappings from a canonical domain onto the physical multiply-connected domain was
proposed for the inverse problem of antiplane strain of a plane with $n$ uniformly stressed inclusions
in \cite{ant2} and \cite{ant3}. The stresses  $\tau_{13}$  and $\tau_{23}$  inside all the inclusions  are equal the same
constants, $\tau_1$ and $\tau_2$, respectively, and are independent of the stresses prescribed at infinity. In \cite{ant2}, the canonical domain was chosen to be a parametric plane with
$n$ slits lying in the real axis, while in \cite{ant3},  the canonical domain was a plane with $n$ circular holes.
In the former case, the method of the Riemnn-Hilbert problem on a hyperelliptic surface was applied
and the conformal map was reconstructed in terms of singular integrals. In the case of a circular domain,
the method of the Riemann-Hilbert problem of the theory of automorphic functions generated by a Schottky symmetry group
was used. 
Since for any doubly and triply connected domain $D_0$ there exists a conformal map 
from a parametric slit domain $\CD$ with slits lying in the same line to the domain $D_0$, 
the method of the slit conformal maps  
was able to recover the whole family
of inclusions by quadratures when $n=2$ or $n=3$ and only a particular family in the case $n\ge 4.$
The circular map gave a series representation of the solution for any finite number of inclusions. 

The goal of this paper is to develop a method of conformal mappings for the inverse antiplane problem with an inclusion
uniformly stressed and  embedded into
a half-plane.  
It is aimed to construct a conformal map that is capable to recover the shape
of the inclusion and do not change the straight boundary of the surrounding semi-infinity body.
 We state the problem and reduce it a boundary value problem for a single analytic function in Section 2.
 In Section 3, we show that if both stresses, $\tau_1=\tau_{13}$ and $\tau_2=\tau_{23}$, are nonzero constants  inside the inclusion,
 then the problem is equivalent to a vector Riemann-Hilbert problem on two contours on an elliptic surface.
 When $\tau_2=0$, the vector problem is decoupled, and two scalar Riemann-Hilbert problems on a finite
 and a semi-infinite contour 
 on the Riemann surface need to be solved. 
 This case is considered in Section 4. To solve these problems, we propose a new analogue of the Cauchy kernel for an elliptic
 surface applicable when the density is not decaying at infinity, while the infinite point lies on the contour of the problem.  
 In Section 5, we  derive the conformal map needed through the solution to  the two Riemann-Hilbert problems solved in the previous
 section and present the results   of our numerical tests.  We generalize the model to the case of $n$ inclusions
 in a half-plane in Section 6. We  show that in two cases of loading the model reduces to vector Riemann-Hilbert problems
 on a  genus-$n$ hyperelliptic surface. We decouple them in the case when the stress $\tau_{23}$ vanishes inside
 all the inclusions.

\setcounter{equation}{0}

\section{Setting}\label{s2}

Consider a semi-infinite elastic body ${\Bbb R}_+^3=\{|x_1|<\infty,  x_2>0, |x_3|<\infty\}$ and an elastic inclusion 
$\{(x_1, x_2)\in D_1, |x_3|<\infty\}$ imbedded into the body. Denote a cross-section of the body
orthogonal to the axis $x_3$ and external to the inclusion by $D_0={\Bbb R}_+^2\setminus D_1$.  
The shear moduli of the domains $D_0$ and $D_1$ are $\mu_0$ and $\mu_1$, respectively, and the inclusion is in ideal contact
with the external body.  
Suppose  the body is subjected to antiplane shear  $\tau_{13}=\tau_1^\infty$ as $x_1\to\pm\infty$, $0<x_2<\infty$.
The boundary of the body $L_0=\{|x_1|<\infty, x_2=0\}$  is free of traction, $\tau_{23}=0$,  and at infinity
as $x_2\to\infty$ and $|x_1|<\infty$, $\tau_{23}=0$. We aim to recover the whole family of possible 
uniformly stressed inclusions $D_1$ such that $\tau_{13}=\tau_1$ and $\tau_{23}=\tau_2$,
$(x_1,x_2)\in D_1$, where $\tau_1$ and $\tau_2$ are prescribed constants.

Denote by $w_j(x_1,x_2)$, $(x_1,x_2)\in D_j$ ($j=0,1$), the $x_3$-component of the displacement vector. The function
$w_j(x_1,x_2)$ is harmonic in $D_j$, and the shear  stresses are expressed through the displacement
by
\beq
\Gs_{13}=\mu_j\fr{\Md w_j}{\Md x_1}, \quad \Gs_{23}=\mu_j\fr{\Md w_j}{\Md x_2}, \quad (x_1,x_2)\in D_j, \quad j=0,1.
\label{2.1}
\eeq
On the interface $L_1$, the boundary conditions of ideal contact read
\beq
w_0=w_1, \quad \mu_0\fr{\Md w_0}{\Md \nu}=\mu_1\fr{\Md w_1}{\Md \nu}, \quad x\in L_1,
\label{2.2}
\eeq
where $\fr{\Md}{\Md \nu}$ is the normal derivative.
Since
\beq 
\mu_1\fr{\Md w_1}{\Md x_1}=\tau_1, \quad \mu_1\fr{\Md w_1}{\Md x_2}=\tau_2,\quad (x_1,x_2)\in D_1,
\label{2.2'}
\eeq
 we can recover the function $w_1$ up to an arbitrary constant $a_1$,
 \beq
 w_1=\fr{\tau_1 x_1+\tau_2 x_2}{\mu_1}+a_1, \quad (x_1,x_2)\in D_1.
 \label{2.3}
 \eeq
 Introduce next two functions, a harmonic conjugate $w_0^*(x_1,x_2)$ of the function $w_0$ and a
 function 
 \beq
 f(z)=w_0(x_1,x_2)+iw_0^*(x_1,x_2)-\fr{\bar\tau z}{\mu_0}, \quad z=x_1+ix_2\in D_0,
 \label{2.4}
 \eeq
 analytic in the domain $D_0$.  Here, $\bar\tau=\tau_1-i\tau_2$. 
 Then, as it was shown in \cite{ant3}, the two real boundary conditions
 (\ref{2.2}) are equivalent  to the following one complex condition on the contour $L_1$ for the analytic function $f(z)$:
 \beq
 f(z)=\fr{1}{\Gl}\R\fr{\bar\tau z}{\mu_0} +a_1+ib_1, \quad z\in L_1,
 \label{2.5}
 \eeq
 where $\Gl=\Gk/(1-\Gk)$, $\Gk=\mu_1/\mu_0$, and $b_1$ is an arbitrary real constant.
By using the asymptotics of $w_0(z)$ as $z\to\infty$ and the Cauchy-Riemann conditions we find
from (\ref{2.4})  that
\beq
f(z)\sim \fr{\tau_1^\infty-\bar\tau}{\mu_0}z+\const, \quad z\to\infty.
\label{2.6}
\eeq
Note that if the  contour $L_1$ had been prescribed, then the boundary condition would constitute an ill-posed problem
since (\ref{2.5}) specifies both the real and the imaginary parts of an analytic function on the contour. 
In the next section we shall apply the method of conformal mappings to recover the function $f(z)$ and the contour $L_1$.
We  shall show that the problem of determination of  the conformal mapping  is equivalent to two Riemann-Hilbert problems
on a genus-1 two-sheeted Riemann surface.

\setcounter{equation}{0}

\section{Vector Riemann-Hilbert problem on an elliptic surface}\label{s3}

Let $z=\Go(\Gz)$ be a conformal map $\Go:\CD\to D_0$ from a parametric $\Gz$-plane cut along
two segments $l_1=[0,1]$ and $l_0=[m,\infty)$ ($m>1$) onto the elastic domain $D_0$.
The function $\Go(\Gz)$ maps the slit $l_1$ sides onto the inclusion boundary $L_1$ and the sides 
of the semi-infinite slit $l_0$ onto the boundary $L_0$ of the elastic half-plane.  

Introduce a new function $F(\Gz)=f(\Go(\Gz))$ analytic in the domain $\CD$. Then the complex boundary condition
(\ref{2.5}) can be equivalently  written in the form
$$
\I F(\Gz)=b_1, \quad \Gz\in l_1,
$$
\beq
\R F(\Gz)=\fr{1}{\Gl}\R\fr{\bar\tau \Go(\Gz)}{\mu_0}+a_1, \quad \Gz\in l_1.
\label{3.1}
\eeq
Derive next  two boundary conditions on the two-sided contour $l_0$. The first condition is obvious. Since $x_2=0$
on the contour $L_0$, we immediately have 
\beq
\I \Go(\Gz)=0, \quad \Gz\in l_0.
\label{3.2}
\eeq
To determine the second relation, we use the condition $\Md w/\Md x_2\to 0$ as $x_2\to 0^+$, $-\infty<x_1<\infty$,    and
the Cauchy-Riemann condition connected the partial derivatives of the function $w$ and its harmonic conjugate $w^*$.
This results in $w^*=b_0$ on the line $L_0$, where $b_0$ is an arbitrary real constant. Therefore
\beq
\I f(z)=b_0-\I\fr{\bar\tau z}{\mu_0}, \quad z\in L_0,
\label{3.3}
\eeq
and the second condition on the slit $l_0$ reads
\beq
\I F(\Gz)=b_0-\I\fr{\bar\tau \Go(\Gz)}{\mu_0}, \quad \Gz\in l_0.
\label{3.4}
\eeq
The boundary conditions (\ref{3.1}), (\ref{3.2}), and (\ref{3.4}) have to be complemented by the conditions at infinity.
We have
\beq
F(\Gz)\sim\fr{\tau_1^\infty-\bar\tau}{\mu_0}\Go(\Gz), \quad \Go(\Gz)\sim c_\pm\Gx^{1/2}, \quad \Gz=\Gx\pm i0, \quad 
\Gx\to\infty,
\label{3.5}
\eeq
where $c_\pm$ are constants.

We wish now to show that  the boundary conditions derived can be transformed into a vector Riemann-Hilbert 
on a Riemann surface.
Consider the algebraic function 
\beq
u^2=p(\Gz), \quad p(\Gz)=\Gz(1-\Gz)(\Gz-m).
\label{3.5'}
\eeq
Fix a single branch of the function $p^{1/2}(\Gz)$ in the $\Gz$-plane cut along the contours $l_0$ and $l_1$
by the condition $p^{1/2}(\Gx)>0$, $-\infty<\Gx<0$. Then on the ``+" and ``-" sides of the cuts $l_0$ and $l_1$,
$\Gz=\Gx\pm i0$ and 
\beq
p^{1/2}(\Gz)=\mp i\sqrt{|p(\Gx)|}, \quad 0<\Gx<1, \quad 
p^{1/2}(\Gz)=\pm i\sqrt{|p(\Gx)|}, \quad m<\Gx<\infty.
\label{3.6}
\eeq
Take two replicas $\CD^+$ and $\CD^-$ of the parametric domain $\CD$ and attach the  ``+" sides of the slits $l_0$
and $l_1$ on the upper sheet $\CD^+$ to the  ``-" sides of the corresponding slits on the lower sheet $\CD^-$.
Then we attach the sides $l^-_j\subset \CD^+$ to the sides $l^+_j\subset \CD^-$, $j=0,1$. 
Points $(\Gz,u(\Gz))$ of the upper sheet and 
$(\bar\Gz,-u(\bar\Gz))$ of the lower sheet  of the resulting genus-1 Riemann surface $\CR$ are symmetric with respect
to the contours $l_0$ and $l_1$. Denote $(\bar\Gz,-u(\bar\Gz)=(\Gz_*,u_*)$. Assume that a point $(\Gz,p^{1/2}(\Gz)$
lies in the upper sheet $\CD^+$. Then the point $(\Gz_*,u_*)$, symmetric to it, lies in the lower sheet $\CD^-$.
On this surface  we introduce two functions
$$
\GF_1(\Gz,u)=\left\{\begin{array}{cc}
F(\Gz), & (\Gz,u)\in\CD^+\subset\CR,\\
\ov{F(\bar\Gz)}, & (\Gz,u)\in\CD^-\subset\CR,\\
\end{array}
\right.
$$
\beq
\GF_2(\Gz,u)=\left\{\begin{array}{cc}
i(\Gl\mu_0)^{-1}\bar\tau\Go(\Gz), & (\Gz,u)\in\CD^+\subset\CR,\\
-i(\Gl\mu_0)^{-1}\tau\ov{\Go(\bar\Gz)}, & (\Gz,u)\in\CD^-\subset\CR.\\
\end{array}
\right.
\label{3.7}
\eeq
These functions are symmetric with the respect to the contour $\CL=l_0\cup l_1\subset \CR$,
\beq
\GF_j(\Gz,u)=\ov{\GF_j(\Gz_*,u_*)}, \quad (\Gx,v)\in\CL,
\label{3.7'}
\eeq
and their limit values on the contour $\CL$ from the upper and lower sheet
are expressed through the functions $F(\Gx)$ and $\Go(\Gx)$ as
$$
\GF_1^+(\Gx,v)=F(\Gx), \quad \GF_2^+(\Gx,v)=i\bar\tau \Go(\Gx),
$$
\beq
\GF_1^-(\Gx,v)=\ov{F(\Gx)}, \quad \GF_2^-(\Gx,v)=-i\tau \ov{\Go(\Gx)},
\label{3.8}
\eeq
where $v=u(\Gx)$.
It remains to write down the boundary conditions (\ref{3.1}), (\ref{3.2}), and (\ref{3.4}) in terms of the functions
(\ref{3.8}).  As a result, we arrive at a vector Riemann-Hilbert problem on the Riemann surface $\CR$ for the  vector
$\GF(\Gz,u)=(\GF_1(\Gz,u), \GF_2(\Gz,u))$. Its boundary condition has the form
\beq
\GF^+(\Gx,v)=G(\Gx,v)\GF^-(\Gx,v)+g(\Gx,v), \quad (\Gx,v)\in\CL\subset\CR,
\label{3.9}
\eeq
where $G(\Gx,v)$ is a piece-wise constant matrix 
\beq
G(\Gx,v)=\left(\begin{array}{cc}
1 & 0\\
2i & 1\\
\end{array}
\right), \quad (\Gx,v)\in l_1, \quad 
G(\Gx,v)=\left(\begin{array}{cc}
1 & i\Gl(1-\bar\tau/\tau)\\
0 & -\bar\tau/\tau\\
\end{array}\right), \quad (\Gx,v)\in l_0,
\label{3.10}
\eeq
and $g(\Gx,v)$ is a piece-wise constant vector
\beq
g(\Gx,v)=\left(\begin{array}{c}
2ib_1\\
-2i(a_1-ib_1) \\
\end{array}
\right), \quad (\Gx,v)\in l_1, \quad 
g(\Gx,v)=\left(\begin{array}{c}
2ib_0 \\
0 \\
\end{array}\right), \quad (\Gx,v)\in l_0.
\label{3.11}
\eeq
The vector $\GF(\Gz,u)$ is symmetric with respect to the contour $\CL$, $\GF(\Gz,u)=\ov{\GF(\Gz_*,u_*)}$, and its components
satisfy the conditions at infinity
\beq
\GF_1^+(\Gz,u)\sim\fr{\Gl(\tau_1^\infty-\bar\tau)}{i\bar\tau}\GF_2^+(\Gz,u), \quad 
\GF_2(\Gz,u)=O(\Gz^{1/2}), \quad \Gz\to\infty.
\label{3.12}
\eeq

\setcounter{equation}{0}

\section{Two scalar Riemann-Hilbert problems on an elliptic surface in the case $\tau_2=0$}\label{s4}

To reduce the vector problem (\ref{3.9}) to two scalar problems to be solved consequently, 
we put $\tau_2=0$. Then $\tau=\bar\tau$,  the matrix $G(\Gx,v)$ is a diagonal matrix on the contour
$l_0$,  and we have two Riemann-Hilbert problems,
\beq
\GF_1^+(\Gx,v)-\GF_1^-(\Gx,v)=2i b_j, \quad (\Gx,v)\in l_j, \quad j=0,1,
\label{4.1}
\eeq
and
$$
\GF_2^+(\Gx,v)+\GF_2^-(\Gx,v)=0, \quad (\Gx,v)\in l_0,
$$
\beq
\GF_2^+(\Gx,v)-\GF_2^-(\Gx,v)=2i[\R\GF_1^+(\Gx,v)-a_1], \quad (\Gx,v)\in l_1.
\label{4.2}
\eeq
The functions $\GF_1(\Gz,u)$ and  $\GF_2(\Gz,u)$ satisfy the symmetry condition (\ref{3.7'}) and the conditions at 
infinity (\ref{3.12}).

\subsection{The first Riemann-Hilbert problem}

For the solution of both Riemann-Hilbert problems, we  propose the following analogue of the Cauchy kernel:
\beq
dV((\Gx,v),(\Gz,u))=\fr12\left[\fr{\Gz-\Gx_0}{\Gx-\Gx_0}+\fr{u}{v}\fr{\Gx-\Gx_0}{\Gz-\Gx_0}\right]
\fr{d\Gx}{\Gx-\Gz}.
\label{4.4}
\eeq
where $u=u(\Gz)$, $v=u(\Gx)$, $\Gx_0$ is an arbitrary real fixed point not lying on the contours $l_0$ and
$l_1$. It will be convenient to take $\Gx_0$ as a negative number.

A singular integral with this kernel on the elliptic surface $\CR$ satisfies the   Sokhotski-Plemelj formulas. The kernel is symmetric
with respect to the contour $\CL$, 
\beq
dV((\Gx,v),(\Gz,u(\Gz)))=\ov{dV((\Gx,v),(\bar\Gz,-u(\bar \Gz)))}.
\label{4.5}
\eeq
If $\Gz$ is a fixed bounded point on the surface $\CR$ and $\Gx\to\pm i0+\infty$, then
$dV$ is decaying as $\Gx^{-3/2}$ that is
\beq
dV\sim \pm\fr{1}{2i}\fr{u(\Gz)}{\Gz-\Gx_0}
|\Gx|^{-3/2}d\Gx, \quad \Gx\to\pm i0+\infty.
\label{4.6}
\eeq
 If $\Gx$ is a fixed point lying either in the contour $l_1$ or any finite part of the contour $l_0$
 and $\Gz\to\infty$, 
 then the kernel $dV$ is bounded,
\beq
dV\sim \left[-\fr{1}{2(\Gx-\Gx_0)}+O(\Gz^{-1/2})\right]d\Gx, \quad \Gz\to \infty.
\label{4.7}
\eeq
  At the two points $(\Gx_0, u(\Gx_0))$ and
 $(\Gx_0, -u(\Gx_0))$ of the surface  the kernel  $dV$ has simple poles. 

Now we shall proceed to solve the Riemann-Hilbert problem (\ref{4.1}).
Owing to the conditions (\ref{3.12}), the general solution 
has the form
$$
\GF_1(\Gz,u)=N_0+i(N_1+iN_2)\fr{u(\Gz)+u(\Gx_0)}{\Gz-\Gx_0}+i(N_1-iN_2)\fr{u(\Gz)-u(\Gx_0)}{\Gz-\Gx_0}
 $$
 \beq
 +\fr{1}{2\pi}\sum_{j=0}^1 b_j \int_{l_j} \left(\fr{\Gz-\Gx_0}{\Gx-\Gx_0}+\fr{u}{v}\fr{\Gx-\Gx_0}{\Gz-\Gx_0}\right)
\fr{d\Gx}{\Gx-\Gz}.
\label{4.8}
\eeq
Here, $N_0$,  $N_1$, and $N_2$ are real constants. At this stage, they are arbitrary. It is directly verified that the 
function (\ref{4.8}) satisfies the symmetry condition (\ref{3.7'}) and at infinity,
\beq
\GF_1(\Gz,u)\sim 2iN_1\fr{u(\Gz)}{\Gz}=O(\Gz^{1/2}), \quad \Gz\to\infty.
\label{4.9}
\eeq
Because of the poles at the points $(\Gx_0,u(\Gx_0))$ and $(\Gx_0,-u(\Gx_0))$ of the kernel $dV$, the solution
(\ref{4.8}) has inadmissible simple poles at these points. These poles become removable singularities if and only if
\beq
2i(N_1+iN_2)+\fr{b_0}{2\pi}\int_{l_0}\fr{d\Gx}{v}+\fr{b_1}{2\pi}\int_{l_1}\fr{d\Gx}{v}=0.
\label{4.10}
\eeq
Computing the two integrals in (\ref{4.10}) we obtain
$$
\int_{l_0}\fr{d\Gx}{v}=\fr{2}{i}\int_m^\infty\fr{d\Gx}{\sqrt{|p(\Gx)|}}=-4ik{\Bbb K}, 
$$
\beq
\int_{l_1}\fr{d\Gx}{v}=-\fr{2}{i}\int_0^1\fr{d\Gx}{\sqrt{|p(\Gx)|}}=
4ik{\Bbb K}, \quad k=m^{-1/2},
\label{4.11}
\eeq
where ${\Bbb K}={\Bbb K}(k)$ is the complete elliptic integral of the first kind. On substituting the values of the integrals into the 
complex equation (\ref{4.10}) we determine
\beq
N_2=0, \quad b_1=b_0-\fr{\pi N_1}{k{\Bbb K}}.
\label{4.12}
\eeq
By removing the zero constant $N_2$ from (\ref{4.8}) we simplify formula (\ref{4.8}),
\beq
\GF_1(\Gz,u)=N_0+\fr{2iN_1u(\Gz)}{\Gz-\Gx_0}
 +\fr{1}{2\pi}\sum_{j=0}^1 b_j \int_{l_j} \left(\fr{\Gz-\Gx_0}{\Gx-\Gx_0}+\fr{u}{v}\fr{\Gx-\Gx_0}{\Gz-\Gx_0}\right)
\fr{d\Gx}{\Gx-\Gz}.
\label{4.13}
\eeq
The function $\GF_1(\Gz,u)$ has three real arbitrary constants, $N_0$, $N_1$, and $b_0$. The fourth constant $b_1$ is expressed through $b_0$ and $N_1$ by equation (\ref{4.12}).

\subsection{The second Riemann-Hilbert problem}

Rewrite equations (\ref{4.2}) in the form
\beq
\GF^+_2(\Gx,v)=G_2(\Gx,v)\GF_2^-(\Gx,v)+g_2(\Gx,v) \quad (\Gx,v)\in\CL,
\label{4.14}
\eeq
where 
\beq
G_2(\Gx,v)=\left\{\begin{array}{cc}
-1, & (\Gx,v)\in l_0,\\
1, & (\Gx,v)\in l_1,\\
\end{array}
\right. \quad
g_2(\Gx,v)=\left\{\begin{array}{cc}
0, & (\Gx,v)\in l_0,\\
2i[\R\GF^+_1(\Gz,v)-a_1], & (\Gx,v)\in l_1.\\
\end{array}
\right.  
\label{4.15}
\eeq
The function $\GF_2(\Gz,u)$ has to satisfy the symmetry condition (\ref{3.7'}) and the condition at infinity (\ref{3.12}).

First we factorize the coefficient $G_2(\Gx,v)$ that is find a function $X(\Gz,u)$ that meets the  symmetry condition $X(\Gz,u)=\ov{X(\Gz_*,u_*)}$, 
$(\Gz,u)\in\CR$,  meromorphic in $\CR\setminus\CL$, and its limit values
satisfy the relation
\beq
G_2(\Gx,v)=X^+(\Gx,v)[X^-(\Gx,v)]^{-1},  \quad (\Gx,v)\in\CL,
\label{4.16}
\eeq
To find such a function, we cut the surface $\CR$ 
along canonical cross-sections $\bf a$ and $\bf b$ to form a simply connected domain $\hat\CR$.
Choose $\bf a$ as a two-sided slit $l_0$ that belongs to both sheets of the surface. The contour  $\bf b$
consists of the segment joining the point $m$ and 1 along the upper sheet $\CD^+$ and the segment $[m,1]\subset\CD^-$.
The positive direction on the boundary of the surface $\hat\CR$, $\Md \hat\CR=\bf a^+ b^+ a^- b^-$, is chosen in the standard way,
that is when a point traverses the boundary, the surface $\hat\CR$ is on the left.    

Show that the function
\beq
X(\Gz,u)=\Gc(\Gz,u)\ov{\Gc(\Gz_*,u_*)} \exp\left\{\fr{1}{2\pi i}\int_{l_0}\log(-1)dV((\Gx,v),(\Gz,u))\right\}
\label{4.17}
\eeq
 solves the factorization problem. Here,
 \beq
 \Gc(\Gz,u)=\exp\left\{-\int_{\Gg} dV((\Gx,v),(\Gz,u))-n_a\int_{\bf a}dV((\Gx,v),(\Gz,u))-n_b\int_{\bf b}dV((\Gx,v),(\Gz,u))\right\},
 \label{4.18}
 \eeq
$n_a$ and $n_b$ are integers, $\Gg$ is a contour lying in the surface $\hat\CR$ and not crossing the canonical cross-sections.
Choose the starting point $q_0=(\Gz_0, u_0)$ of the contour $\Gg$ as a point on the upper
sheet $\CD^+$,  $u_0=\sqrt{p(\Gz_0)}$. 
The terminal point $q_1=(\Gz_1, u_1)$, $u_1=u(\Gz_1)$,  cannot be chosen a priori. It may fall on either sheet
 and is to be determined.

By the Sokhotski-Plemelj formulas, the function $X(\Gz,u)$ satisfies the relation (\ref{4.16}).
We also assert that the function $X(\Gz,u)$ is symmetric with respect to the contour $\CL$, $X(\Gz,u)=\ov{X(\Gz_*,u_*)}$,
$(\Gz,u)\in\CR$. Although the integrals in (\ref{4.18}) have jumps multiple of $2\pi i$ across the contours
of integration, the functions $\Gc(\Gz,u)$ and  $\Gc(\Gz_*,u_*)$ are meromorphic on the surface $\CR$.
Owing to logarithmic singularities of the integral over the contour $\Gg$ in (\ref{4.18}) at the endpoints of the contour, the function
$\Gc(\Gz,u)$ has a simple zero at the point $q_0\in \CD^+$ and a simple pole at the point $q_1\in\CR$.
Because of the symmetry, the factorizing function $X(\Gz,u)$  has  two simple zeros at the points $q_0$ and 
$q_{0*}=(\bar\Gz_0, -\sqrt{p(\bar\Gz_0)})$ and two simple poles 
 at   the points $q_1$ and 
$q_{1*}=(\bar\Gz_1, -u(\bar\Gz_1))$. 

In view of the two simple poles at the points $(\Gx_0, u(\Gx_0))$ and
 $(\Gx_0, -u(\Gx_0))$ of the kernel  $dV$,  the function $X(\Gz,u)$ has essential singularities  at these points.
 These points become removable singular points if the following condition is satisfied: 
\beq
i\left(\fr14-n_a\right)\int_{l_0}\fr{d\Gx}{u(\Gx)}-\fr{i}{2}\int_\Gg\fr{d\Gx}{u(\Gx)}+\fr{i}{2}\int_\Gg\fr{d\bar\Gx}{u(\bar\Gx)}=0.
\label{4.19}
\eeq
This relation is the imaginary equation of the Jacobi inversion problem on the elliptic surface $\CR$ with respect
to the point $q_1=(\Gz_1,u_1)\in\CR$ and the integers $n_a$ and $n_b$ \cite{az}
\beq
\int_0^{q_1}\fr{d\Gx}{u(\Gx)}+n_a A+n_b B=h,
\label{4.20}
\eeq
where $A$ and $B$ are the cyclic periods of the abelian integral in (\ref{4.20})
$$
A=\int_{\bf a}\fr{d\Gx}{u(\Gx)}=-2i\int_m^\infty\fr{d\Gx}{\sqrt{|p(\Gx)|}}=-4ik{\Bbb K},
$$
\beq
B=\int_{\bf b}\fr{d\Gx}{u(\Gx)}=2\int_1^m\fr{d\Gx}{\sqrt{|p(\Gx)|}}=4k{\Bbb K}',
\label{4.21}
\eeq
${\Bbb K}'={\Bbb K}(\sqrt{1-k^2})$, and  $h$ is given by
 \beq
 h=\int_0^{\Gz_0}\fr{d\Gx}{p^{1/2}(\Gx)}-ik{\Bbb K}.
 \label{4.22}
 \eeq
The affix of the point $q_1$ is determined in terms of the elliptic sine
\beq
\Gz_1={\rm sn}^2\fr{ih}{2k}.
\label{4.23}
\eeq
To write down the values of the integers $n_a$ and $n_b$, introduce two constants
\beq
\CI_\pm=\int_0^{\Gz_0}\fr{d\Gx}{p^{1/2}(\Gx)}\pm \int_0^{\Gz_1}\fr{d\Gx}{p^{1/2}(\Gx)}-ik {\Bbb K}.
\label{4.24}
\eeq
If both constants computed by
\beq
n_a=-\fr{\I \CI_-}{4k{\Bbb K}}, \quad n_b=\fr{\R\CI_-}{4k{\Bbb K'}}
\label{4.25}
\eeq
are integers, then the point $q_1=(\Gz_1,\sqrt{p(\Gz_1)})\in\CD^+$. Otherwise, if at least one of the right-hand 
sides in equations (\ref{4.25}) is not integer, then the point
 $q_1=(\Gz_1,-\sqrt{p(\Gz_1)})$ falls on the lower sheet $\CD^-$, and  the constants $n_a$ and $n_b$ defined by
 \beq
n_a=-\fr{\I \CI_+}{4k{\Bbb K}}, \quad n_b=\fr{\R\CI_+}{4k{\Bbb K'}}
\label{4.25'}
\eeq
have to be integers.

It turns out that owing to the symmetry of the Riemann surface, the factorizing function $X(\Gz,u)$
is independent of the integer $n_b$, and formula (\ref{4.17}) can be simplified and written in the form
$$
X(\Gz,u)=\exp
\left\{\left(\fr12-2n_a\right)\fr{u(\Gz)}{i(\Gz-\Gx_0)}
\int_m^\infty \fr{(\Gx-\Gx_0)d\Gx}{\sqrt{|p(\Gx)|}(\Gx-\Gz)}
\right.
$$
\beq
\left.
-\fr12\int_\Gg\left(\fr{\Gz-\Gx_0}{\Gx-\Gx_0}+\fr{u(\Gz)}{u(\Gx)}\fr{\Gx-\Gx_0}{\Gz-\Gx_0}\right)
\fr{d\Gx}{\Gx-\Gz}
-\fr12\int_\Gg\left(\fr{\Gz-\Gx_0}{\bar\Gx-\Gx_0}-\fr{u(\Gz)}{u(\bar\Gx)}\fr{\bar\Gx-\Gx_0}{\Gz-\Gx_0}\right)\fr{d\bar\Gx}{\bar\Gx-\Gz}
\right\}.
\label{4.26}
\eeq

Having factorized the function $G_2(\Gx,v)$ we can proceed with the solution of the Riemann-Hilbert problem (\ref{4.14}).
On replacing the function $G_2(\Gx,v)$ in (\ref{4.14}) by $X^+(\Gx,v)[X^-(\Gx,v)]^{-1}$ we arrive at
\beq
\fr{\GF^+_2(\Gx,v)}{X^+(\Gx,v)}=\fr{\GF^-_2(\Gx,v)}{X^-(\Gx,v)}+\fr{g_2(\Gx,v)}{X^+(\Gx,v)}, \quad (\Gx,v)\in\CL.
\label{4.27}
\eeq
Introduce next a singular integral
\beq
\Psi(\Gz,u)=\fr{1}{2\pi i}\int_{l_1}\fr{g_2(\Gx,v)}{X^+(\Gx,v)}dV((\Gx,v),(\Gz,u)), \quad (\Gz,u)\in\CR.
\label{4.28}
\eeq
It will be convenient to represent this integral over the two-sided contour $l_1$ as follows:
$$
\Psi(\Gz,u)=\fr{1}{4\pi i}\int_0^1\left[\fr{g_2(\Gx^+,v^+)}{X^+(\Gx^+,v^+)}
\left(
\fr{\Gz-\Gx_0}{\Gx-\Gx_0}+\fr{u}{v^+}\fr{\Gx-\Gx_0}{\Gz-\Gx_0}\right)
\right.
$$
\beq
\left.
-\fr{g_2(\Gx^-,v^-)}{X^+(\Gx^-,v^-)}
\left(
\fr{\Gz-\Gx_0}{\Gx-\Gx_0}+\fr{u}{v^-}\fr{\Gx-\Gx_0}{\Gz-\Gx_0}\right)
\right]
\fr{d\Gx}{\Gx-\Gz}.
\label{4.29}
\eeq
Here, $(\Gx^\pm,v^\pm)=(\Gx\pm i0, \sqrt{p(\Gx\pm i0)}$ are points on the sides $l_1^\pm\subset\CD^+$, and $u=u(\Gz)$.
The values $g_2(\Gx^\pm,v^\pm)$ are recovered from (\ref{4.13}) and (\ref{4.15}). We have
\beq
g_2(\Gx^\pm,v^\pm)=2i\left[N_0^*\pm\fr{\sqrt{|p(\Gx)|}}{\Gx-\Gx_0}g_0(\Gx)\right], \quad 0<\Gx<1,
\label{4.30}
\eeq
where $N_0^*=N_0-a_1$ and
\beq
g_0(\Gx)=2N_1-\fr{b_0}{\pi}\int_m^\infty\fr{(\tau-\Gx_0)d\tau}{\sqrt{|p(\tau)|}(\tau-\Gx)}+
\fr{b_1}{\pi}\int_0^1\fr{(\tau-\Gx_0)d\tau}{\sqrt{|p(\tau)|}(\tau-\Gx)}, \quad 0<\Gx<1.
\label{4.31}
\eeq
The Cauchy principal value for the last singular integral is assigned.

We now use the Sokhotski-Plemelj formulas to represent the second term in the right-hand side of equation 
(\ref{4.27}) by the difference of the limit values of the function $\Psi(\Gz,u)$, $\Psi^+(\Gx,v)-\Psi^-(\Gx,v)$,
apply the continuity principle and the generalized Liouville theorem on the surface $\CR$ and deduce
the general solution to the Riemann-Hilbert problem (\ref{4.14})
\beq
\GF_2(\Gz,u)=X(\Gz,u)[\Psi(\Gz,u)+\GO(\Gz,u)], \quad (\Gz,u)\in\CR.
\label{4.32}
\eeq
Here, $\GO(\Gz,u)$ is a rational function of the surface $\CR$. It has to have the following properties.

(i) In virtue of the symmetry of the functions $\GF_2(\Gz,u)$, $X(\Gz,u)$, and $\Psi(\Gz,u)$
with respect to the contour $\CL$, the function $\GO(\Gz,u)$ has to be symmetric with respect to this contour as well.

(ii) Owing to the simple pole of the kernel $dV$ at $\Gz=\Gx_0$ on both sheets of the surface $\CR$, the 
function $\GO(\Gz,u)$ also has to have the poles at these points and
\beq
 \mathop{\rm res}\limits_{\Gz=\Gx_0}[\Psi(\Gz,u)+\GO(\Gz,u)]=0.
 \label{4.33}
 \eeq
 
 (iii) Since the factorizing function $X(\Gz,u)$ has a simple pole at the point $q_1=(\Gz_1,u_1)$ that
 lies either on the upper or lower sheet of the surface (this is determined by the solution of the Jacobi problem),
 \beq
\Psi(\Gz_1,u_1)+\GO(\Gz_1,u_1)]=0.
 \label{4.34}
 \eeq

(iv)  Because of the simple zero of the function $X(\Gz,u)$ at the point $q_0=(\Gz_0,u_0)\in\CD^+$, 
the function $\GO(\Gz,u)$ may have a simple pole at this point.

(v) $\GO(\Gz,u)=O(\Gz^{1/2})$, $\Gz\to\infty$, and the principal term of its asymptotics  at infinity
has to be chosen such that the first condition in (\ref{3.12}) holds.

The most general form of the rational function $\GO(\Gz,u)$  has the form
\beq
\GO(\Gz,u)=M_0+(M_1+iM_2)\fr{u(\Gz)+u(\Gz_0)}{\Gz-\Gz_0}-
(M_1-iM_2)\fr{u(\Gz)-u(\bar\Gz_0)}{\Gz-\bar\Gz_0}+\fr{2iM_3u(\Gz)}{\Gz-\Gx_0},
\label{4.35}
\eeq
where $M_j$ ($j=0,\ldots,3$) are arbitrary real constants.  The function (\ref{4.35}) is symmetric with respect
to the contour $\CL$ and has a simple pole at the point $q_0=(\Gz_0,u_0)\in\CD^+$. Thus, properties (i) and (iv)
have been satisfied. 

To meet the condition (v), that is the first relation in (\ref{3.12}) at infinity, we compute the limit
\beq
\lim_{\Gz\to\infty} X(\Gz,u)=iX_\infty, \quad X_\infty=\left|\fr{\Gz_1-\Gx_0}{\Gz_0-\Gx_0}\right|
\label{4.36}
\eeq
and pass to the limit $\Gz\to\infty$ in formulas (\ref{4.13}) and (\ref{4.32}).  We deduce
\beq
M_2=-M_3+\fr{\tau_1N_1}{\Gl(\tau_1^\infty-\tau_1)X_\infty}.
\label{4.37}
\eeq
We next compute the residue in (\ref{4.33}), and the property (ii)
yields the constant $M_3$
\beq
M_3=-\fr{1}{8\pi i}\int_0^1\left[
\fr{g_2(\Gx^+,v^+)}{X^+(\Gx^+,v^+)}+\fr{g_2(\Gx^-,v^-)}{X^+(\Gx^-,v^-)}
\right]\fr{d\Gx}{\sqrt{|p(\Gx)|}}.
\label{4.38}
\eeq
 The last property to be satisfied is the condition  (iii).
 Denote
 $$
 P+iQ=\Psi(\Gz_1,u_1), \quad P_0+iQ_0=\fr{2u(\Gz_1)}{\Gz_1-\Gx_0},
 $$
 \beq
 P_1+iQ_1=\fr{u(\Gz_1)+u(\Gz_0)}{\Gz_1-\Gz_0},\quad 
  P_2+iQ_2=\fr{u(\Gz_1)-u(\bar\Gz_0)}{\Gz_1-\bar\Gz_0}.
  \label{4.46}
  \eeq
  where $P$, $Q$, $P_l$, and $Q_l$ ($l=0,1,2$) are real constants.
Then equation (\ref{4.34}) implies 
$$
M_1=\fr{(P_1+P_2)M_2+P_0 M_3+Q}{Q_2-Q_1},
$$
\beq
M_0=(P_2-P_1)M_1+(Q_1+Q_2)M_2+Q_0M_3-P.
\label{4.47}
\eeq

\setcounter{equation}{0}

\section{Family of conformal maps. Numerical results}

The family of conformal mappings $z=\Go(\Gz)$ from the doubly connected $\Gz$-domain $\CD$ onto 
the doubly-connected  domain $D_0$, the exterior of the inclusion  $D_1$  embedded into a $z$-plane, is
described by the formula
\beq
\Go(\Gz)=-\fr{i\Gl}{\hat\tau_1}\GF_2^+(\Gz,u), \quad (\Gz, u)\in \CD^+.
\label{5.1}
\eeq
In addition to the three  real parameters of the elasticity model,
$\hat\tau_1=\tau_1/\mu_0$,  $\hat\tau_1^\infty=\tau_1^\infty/\mu_0$, and $\Gk=\mu_1/\mu_0$, the map possesses four other real parameters. They are $N_0^*$, $N_1\ne 0$, $b_0$, and $m\in(1,\infty)$. 
The nonzero real parameter $N_1$ is a scaling parameter, and 
the map $\hat\Go(\Gz)=N_1^{-1}\Go(\Gz)$ has three real parameters, $\hat N_0^*=N^*_0/N_1$,
$\hat b_0=b_0/N_1$,  and the geometric parameter $m$. Our numerical tests show that 
the map is invariant of the  second parameter $\hat b_0$. As for the first parameter $\hat N_0^*$, it is a translation parameter.
Variation of this parameter leads to translation of the inclusion along the $x$-axis and does not change the inclusion profile
and the distancee of the inclusion points to the $x$-axis. 
This means that given the model three parameters
$\hat\tau_1$,  $\hat\tau_1^\infty$, and  $\Gk$, the parameter $m\in(1,\infty)$
generates a one-parametric family of maps $\hat\Go(\Gz)$  and therefore a one-parametric family of scaled 
uniformly stressed inclusions  embedded into a half-plane. The parameter $m>1$ has to be chosen such that 
the inclusion boundary does not intersect the $x$-axis, the boundary of the external elastic body.
 
 To verify that the function $z=\Go(\Gz)$ maps the two-sided slits $l_0$ and $l_1$ onto the boundary of the $z$-half-plane
 and the inclusion $D_1$, respectively,  indeed, we write down the function (\ref{5.1}) on the contours $l_0^\pm$
 and $l_1^\pm$. We have
 \beq
 \Go(\Gx^\pm)=-\fr{i\Gl}{\hat\tau_1}X^+(\Gx^\pm,u^\pm)[\Psi(\Gx^\pm,u^\pm)+\GO(\Gx^\pm,u^\pm)], \quad (\Gx^\pm,u^\pm)\in l_0^\pm\cup l_1^\pm\subset \CD^+,
 \label{5.2}
 \eeq
where $\Gx^\pm=\Gx\pm i0$, $u^\pm=p^{1/2}(\Gx^\pm)$. On the two sides of the contour $l_0^\pm$, the functions $\Psi(\Gx^\pm,u^\pm)$ and $\GO(\Gx^\pm,u^\pm)$ are real. By the Sokhotski-Plemelj
 formulas we discover  that $X^+(\Gx^\pm, u^\pm)=i\R X^+(\Gx^\pm, u^\pm)$,
and therefore the imaginary part of the function $\Go(\Gx^\pm)$ is equal to 0, while $\R\Go(\Gx^\pm)\in(-\infty,\infty)$.
This means that the contour $l_0$ maps to the real axis of the physical plane.

Consider now the contour $l_1$. On applying the Sokhotski-Plemelj formulas  we obtain from (\ref{4.13}), (\ref{4.32}) that
 $$
 \GF_1^+(\Gx,v)-\GF_1^-(\Gx,v)=2i b_1, \quad (\Gx,v)\in l_1,
 $$
 \beq
 \GF_2^+(\Gx,v)-\GF_2^-(\Gx,v)=2i [\R\GF_1^+(\Gx,v)-a_1], \quad (\Gx,v)\in l_1.
 \label{5.3}
 \eeq
In view of (\ref{3.7}) we conclude that  the complex boundary condition (\ref{2.5})
and the interface conditions (\ref{2.2}) hold.
It is directly verified that the solution determined satisfies the condition at infinity (\ref{2.6}).

To recover the contour $L_1$, the boundary of the inclusion, we need to let a point $\Gz$ traverse the contour
$l_1$ along the positive and negative sides. 
Since all our numerical tests show that the point $\Gz_1$ falls on the upper sheet of the surface $\CR$,
we simplify the formula for the function $X^+(\Gx^\pm,u^\pm)$ for this case
$$
X(\Gx^\pm,u^\pm)=\exp
\left\{\mp\left(\fr12-2n_a\right)\fr{\sqrt{|p(\Gx)|}}{\Gx-\Gx_0}
\int_0^{1/m} \fr{(1-\Gj\Gx_0)d\Gj}{\sqrt{\Gj(1-\Gj)(1-m\Gj)}(1-\Gj\Gx)}
\right.
$$
\beq
\left.
-\R\left[
\int_\Gg
\left(\fr{\Gx-\Gx_0}{\Gj-\Gx_0}
\mp\fr{i\sqrt{|p(\Gx)|}}{\sqrt{p(\Gj)}}\fr{\Gj-\Gx_0}{\Gx-\Gx_0}\right)
\fr{d\Gj}{\Gj-\Gx}
\right]
\right\}, \quad (\Gx^\pm,v^\pm)\in l_1^\pm\subset\CD^+.
\label{5.4}
\eeq
A similar formula can be written when $\Gz_1\in\CD^-$. In this case  a part of the contour $\Gg$, a contour $(\Gz_0,0)$,
lies on the upper sheet
$\CD^+$
and the second part $(0,\Gz_1)$ lies on the lower sheet $\CD^-$.
Both of the integrals in (\ref{5.4}) are nonsingular, and the Gauss quadrature formulas give a good accuracy of computations.

Write next the limit values of  the function $\Psi(\Gz,u)$ on the sides $l_1^\pm$ of the contour $l_1$. By the Sokhotski-Plemelj
formulas we find from (\ref{4.29}) 
$$
\Psi(\Gx^\pm,u^\pm)=\fr{g_2(\Gx^\pm,u^\pm)}{2X^+(\Gx^\pm,u^\pm)}+
\fr{1}{4\pi i}\int_0^1\left[\fr{g_2(\Gj^+,v^+)}{X^+(\Gj^+,v^+)}
\left(
\fr{\Gx-\Gx_0}{\Gj-\Gx_0}+\fr{u^\pm}{v^+}\fr{\Gj-\Gx_0}{\Gx-\Gx_0}\right)
\right.
$$
\beq
\left.
-\fr{g_2(\Gj^-,v^-)}{X^+(\Gj^-,v^-)}
\left(
\fr{\Gx-\Gx_0}{\Gj-\Gx_0}+\fr{u^\pm}{v^-}\fr{\Gj-\Gx_0}{\Gx-\Gx_0}
\right)
\right]
\fr{d\Gj}{\Gj-\Gx}, \quad (\Gx^\pm,u^\pm)\in l_1^\pm.
\label{5.5}
\eeq
 On replacing the function $g_2$ by its expression (\ref{4.30}) 
we split the integral into three other integrals. We write
\beq
\Psi(\Gx^\pm,u^\pm)=\fr{i}{X^+(\Gx^\pm,u^\pm)}\left[N_0^*\pm\fr{\sqrt{|p(\Gx)|}}{\Gx-\Gx_0}g_0(\Gx)\right]
+\CI_1(\Gx)\pm \CI_2(\Gx)+\CI_3^\pm(\Gx),\quad 0<\Gx<1,
\label{5.6}
\eeq
where 
$$
\CI_1(\Gx)=\fr{\Gx-\Gx_0}{2\pi}\int_0^1\fr{\sqrt{|p(\Gj)|}g_0(\Gj)}{(\Gj-\Gx_0)^2}
Y_+(\Gj)
\fr{d\Gj}{\Gj-\Gx},
$$
$$
\CI_2(\Gx)=\fr{N_0^*\sqrt{|p(\Gx)|}}{2\pi(\Gx-\Gx_0)}
\int_0^1\fr{\Gj-\Gx_0}{\sqrt{|p(\Gj)|}}Y_+(\Gj)
\fr{d\Gj}{\Gj-\Gx},
$$
\beq
\CI_3^\pm(\Gx)=\fr{1}{2\pi}
\int_0^1\left[
\fr{N_0^*(\Gx-\Gx_0)}{\Gj-\Gx_0}\pm 
\fr{\sqrt{|p(\Gx)|}g_0(\Gj)}{\Gx-\Gx_0}\right]
Y_-(\Gj)
\fr{d\Gj}{\Gj-\Gx},
\label{5.7}
\eeq
and
\beq
Y_\pm(\Gj)=\fr{1}{X^+(\Gj^+,v^+)}\pm\fr{1}{X^+(\Gj^-,v^-)}.
\label{5.8}
\eeq
We assert that all these three integrals are singular. The density of the first integral vanishes at the endpoints,
and the integral can be represented in the form
\beq
\CI_1(\Gx)=\fr{\Gx-\Gx_0}{2\pi}\int_0^1
\fr{\sqrt{\Gj(1-\Gj)}\CF_1(\Gj, \Gx)d\Gj}{\Gj-\Gx},
\label{5.9}
\eeq
where
\beq
\CF_1(\Gj,\Gx)=\fr{\sqrt{m-\Gj}g_0(\Gj)Y_+(\Gj)}{(\Gj-\Gx)^2}.
\label{5.10}
\eeq
To compute this integral, we expand the function $\CF_0$ in terms of the Chebyshev polynomials of the second kind
\beq
\CF_1(\Gj,\Gx)=\sum_{l=1}d_l(\Gx) U_{l-1}(2\Gj-1), \quad 0<\Gx<1,
\label{5.11}
\eeq
where the expansion coefficients 
\beq
d_l(\Gx)=\fr{2}{\pi}\int_{-1}^1 \CF_1\left(\fr{\Gj+1}{2},\Gx\right)U_{l-1}(\Gj)\sqrt{1-\Gj^2}d\Gj
\label{5.12}
\eeq
are evaluated by the order-$N$ Gauss quadrature formula
\beq
d_l(\Gx)=\fr{2}{N+1}\sum_{j=1}^N\sin\fr{j\pi}{N+1}\sin\fr{j l\pi}{N+1}\CF_1\left(\fr{x^0_j+1}{2},\Gx\right),\quad x_j^0=\cos\fr{j\pi}{N+1}.
\label{5.13}
\eeq
On substituting the expansion (\ref{5.11}) into formula (\ref{5.9}) and using the integral relation for the Chebyshev polynomials 
of the first and second kind 
\beq
\int_{-1}^1\fr{\sqrt{1-\Gj^2}U_{l-1}(\Gj)d\Gj}{\Gj-\Gx}=-\pi T_l(\Gx), \quad |\Gx|<1, \quad l=1,2,\ldots,
\label{5.14}
\eeq
we deduce the series expansion for the integral (\ref{5.9}) 
\beq
\CI_1(\Gx)=-\fr{\Gx-\Gx_0}{4}\sum_{l=1}^\infty d_l (\Gx)T_l(2\Gx-1), \quad 0<\Gx<1.
\label{5.15}
\eeq
The density of the second integral $\CI_2(\Gx)$ in (\ref{5.7}) has the square root singularity at the endpoints.
The integral can be computed in a similar fashion \cite{ant2}
\beq
\CI_2(\Gx)= \fr{N_0^*\sqrt{|p(\Gx)|}}{\Gx-\Gx_0}\sum_{l=1}^\infty c_lU_{l-1}(2\Gx-1),
\label{5.16}
\eeq
where
$$
c_l=\fr{2}{N}\sum_{j=1}^N\CF_2\left(\fr{1+x_j}{2}\right)\cos\fr{(2j-1)l\pi}{2N},
$$
\beq
x_j=\cos\fr{(2j-1)\pi}{2N}, \quad 
\CF_2(\Gj)=\fr{(\Gj-\Gx_0)Y_+(\Gj)}{\sqrt{m-\Gj}}.
\label{5.17}
\eeq
Owing to the function $Y_-(\Gj)$, the density of the integral $\CI_3(\Gx)$ in  (\ref{5.7}) vanishes at the endpoints,
and both formulas employed for the integrals $\CI_1(\Gx)$ and $\CI_2(\Gx)$ give a good accuracy.

Note that the first term in formula
(\ref{5.6}) and the integrands of the integrals $\CI_1(\Gx)$ and $\CI_3(\Gx)$ have the function $g_0(\Gj)$
that has two  integrals itself. One of them is not singular, while the second one is understood in the principal value sense,
\beq
g_0(\Gx)=2N_1  -\fr{b_0}{\pi\sqrt{m}}\int_0^{1/m}
\fr{\tilde\CF_0(\Gj,\Gx)d\Gj}{\sqrt{\Gj(1/m-\Gj)}}
+\fr{b_1}{\pi}\int_0^1\fr{\tilde\CF_1(\Gj)d\Gj}{\sqrt{\Gj(1-\Gj)}(\Gj-\Gx)},\quad 0<\Gx<1.
\label{5.18}
\eeq
Here,
\beq
\tilde\CF_0(\Gj,\Gx)=\fr{1-\Gj\Gx_0}{\sqrt{1-\Gj}(1-\Gj\Gx)},
\quad \tilde\CF_1(\Gj)=\fr{\Gj-\Gx_0}{\sqrt{m-\Gj}}.
\label{5.19}
\eeq
The first integral in (\ref{5.18}) is evaluated by the Gauss quadrature rule, while the second one is computed by 
expanding it in terms of the Chebyshev polynomials of the second kind as it was done for the integral $\CI_2(\Gx)$
in (\ref{5.16}).

\begin{figure}[t]
\centerline{
\scalebox{0.6}{\includegraphics{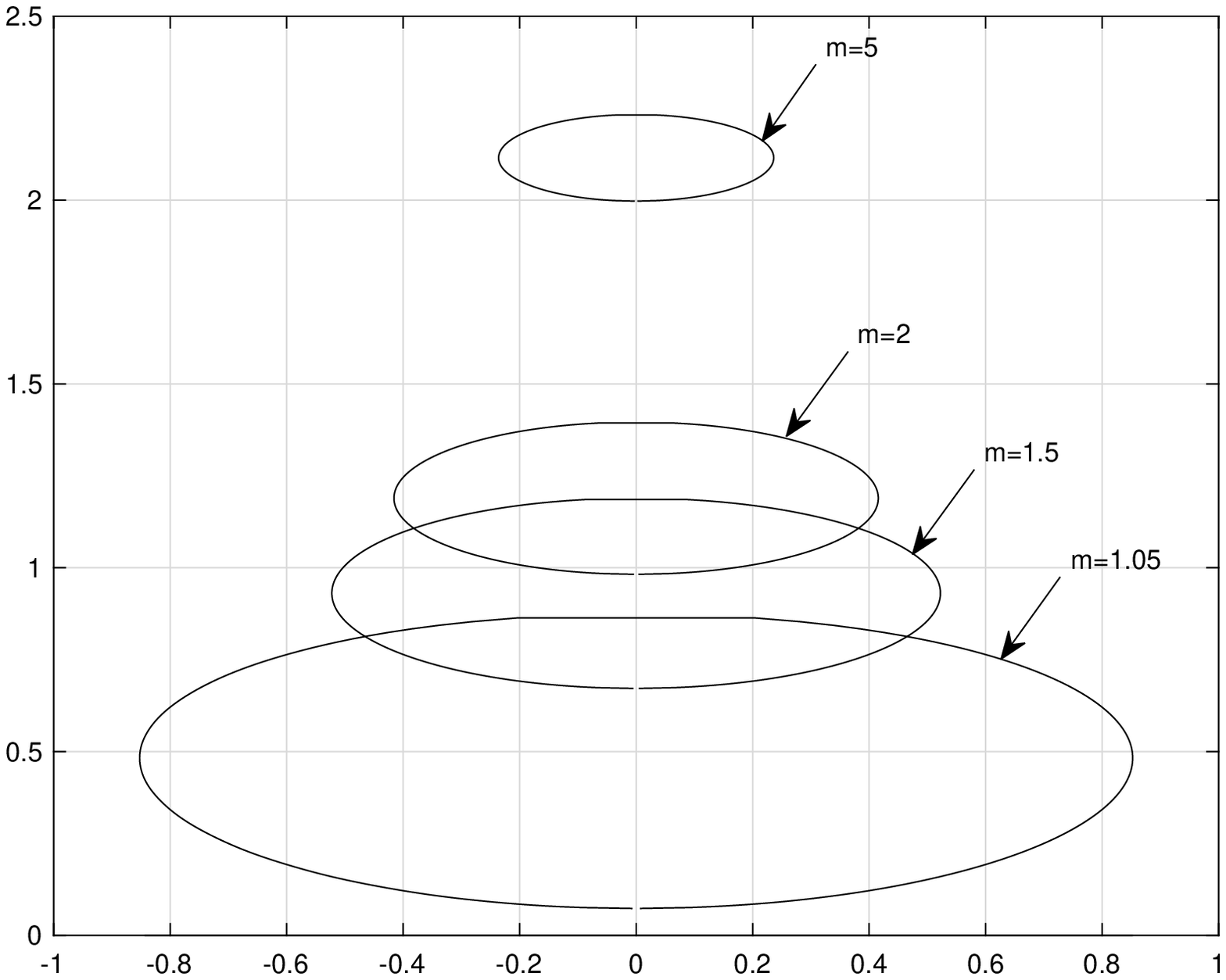}}}
\caption{Normalized inclusion in the half-plane $|x|<\infty$, $y\ge 0$ for different values of $m$ when
$\Gk=0.5$, $N_0^*=0$, $\tau_1/\mu_0=-1$, $\tau_1^\infty/\mu_0=-2$.} 
\label{fig1}
\end{figure} 

\begin{figure}[t]
\centerline{
\scalebox{0.6}{\includegraphics{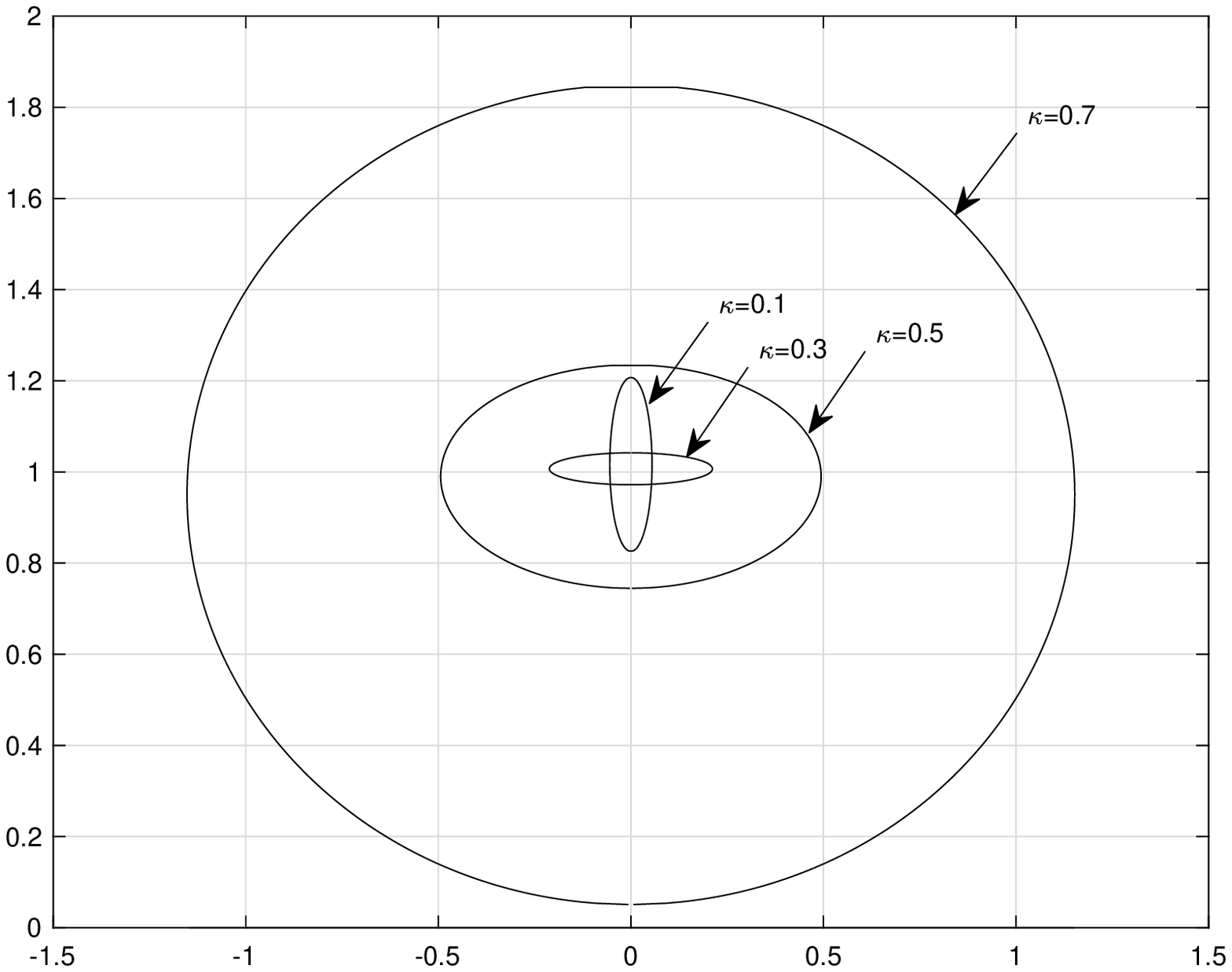}}}
\caption{Normalized inclusion in the half-plane $|x|<\infty$, $y\ge 0$ for different values of $\kappa\in(0,1)$ when
$m=1.6$, $N_0^*=0$, $\tau_1/\mu_0=-1$, $\tau_1^\infty/\mu_0=-2$.} 
\label{fig2}
\end{figure} 

\begin{figure}[t]
\centerline{
\scalebox{0.6}{\includegraphics{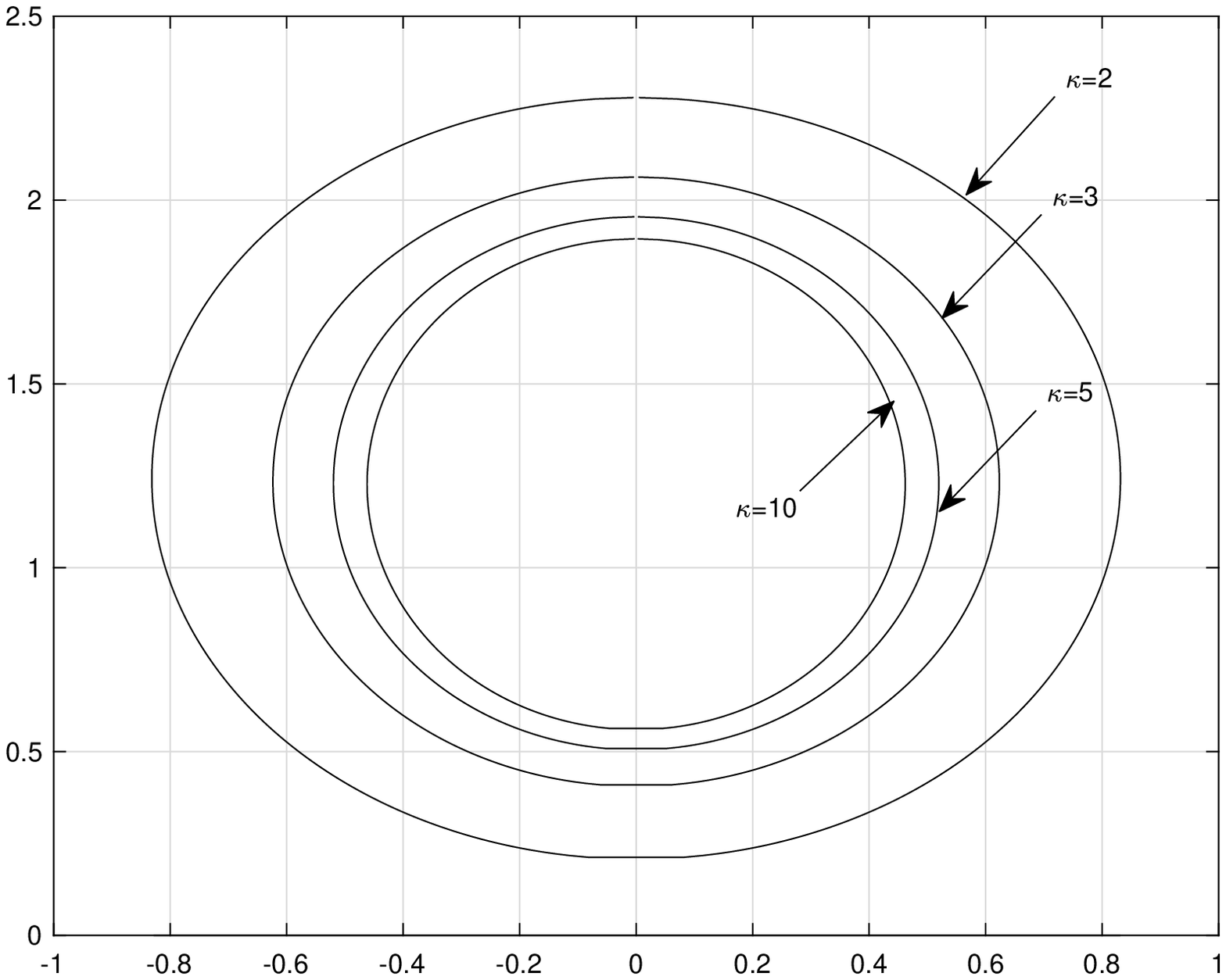}}}
\caption{Normalized inclusion in the half-plane $|x|<\infty$, $y\ge 0$ for different values 
of $\kappa>1$ when
$m=2$, $N_0^*=0$,  $\tau_1/\mu_0=-1$, $\tau_1^\infty/\mu_0=-2$.} 
\label{fig3}
\end{figure}

\begin{figure}[t]
\centerline{
\scalebox{0.6}{\includegraphics{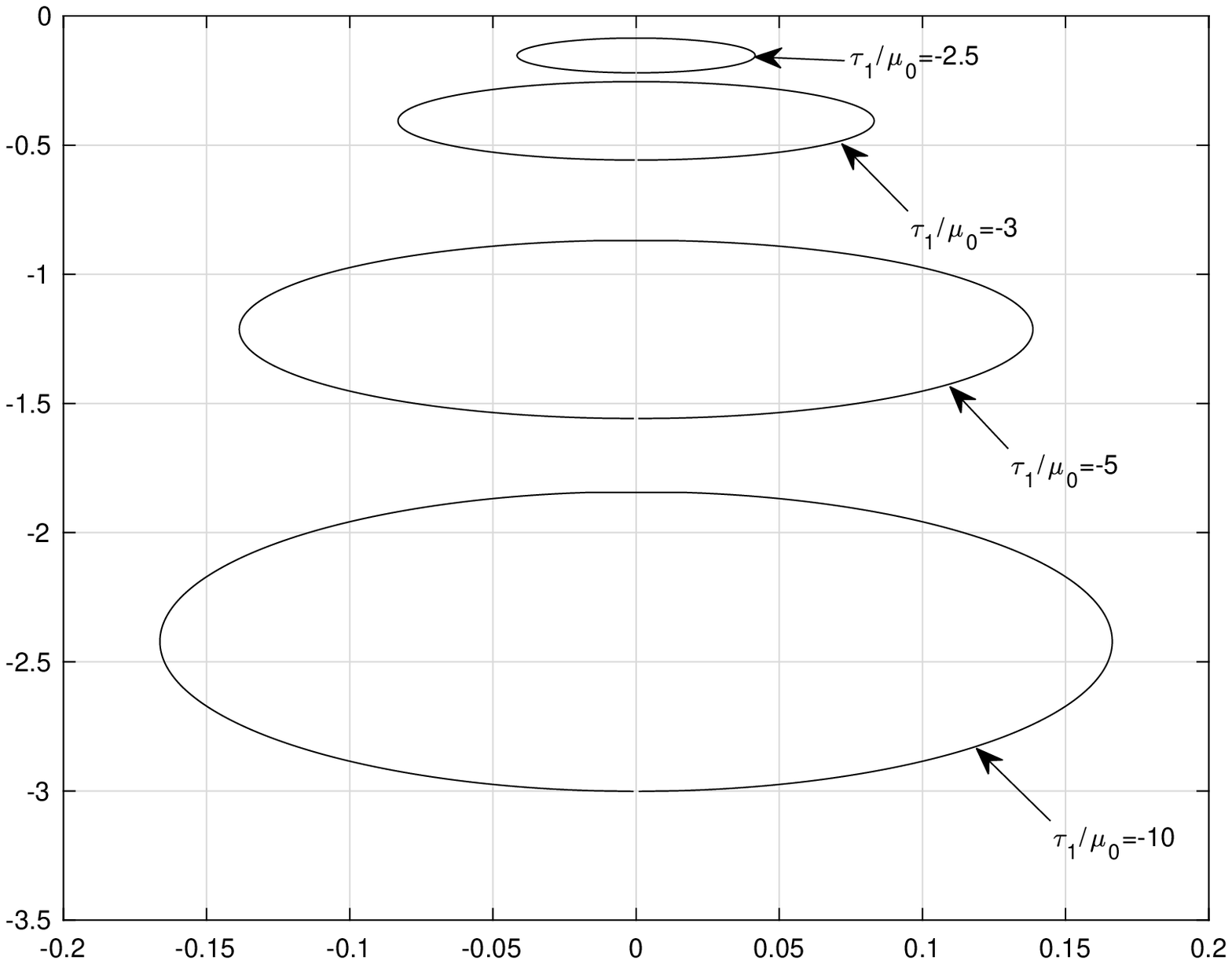}}}
\caption{Normalized inclusion in the half-plane $|x|<\infty$, $y\le 0$ for different values of $\tau_1/\mu_0$ when
$m=2$, $N_0^*=0$, $\Gk=0.5$, $\tau_1^\infty/\mu_0=-2$.} 
\label{fig4}
\end{figure}

\begin{figure}[t]
\centerline{
\scalebox{0.6}{\includegraphics{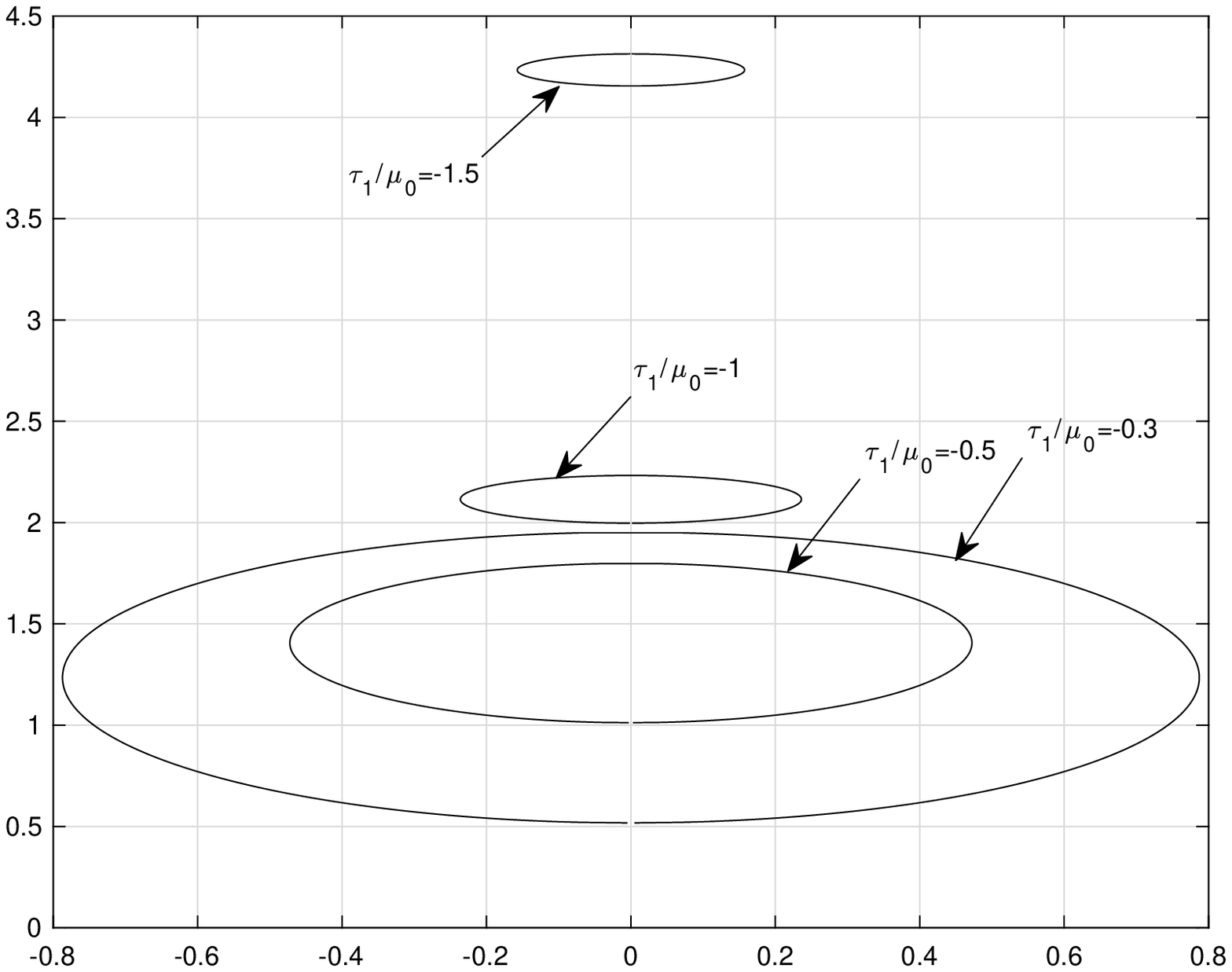}}}
\caption{Normalized inclusion in the half-plane $|x|<\infty$, $y\ge 0$ for different values of $\tau_1/\mu_0$ when
$m=5$, $N_0^*=0$, $\Gk=0.5$, $\tau_1^\infty/\mu_0=-2$.} 
\label{fig5}
\end{figure}

The results of computations are presented in Figures 1 to 5. In Fig. 1, for the model parameters selected to be
$\Gk=\mu_1/\mu_0=0.5$, $\tau_1/\mu_0=-1$, and $\tau_1^\infty/\mu_0=-2$ we plot the inclusion profile
for some values of the mapping parameter $m$. By decreasing the parameter $m$ we increase the inclusion size
and eventually intersect the boundary of the external body. At the same time, given a set of the problem parameters
$\Gk\ne 1$ and $\tau_1\ne\tau_1^\infty$, it is always possible to increase the parameter $m$ such that the inclusion
is completely embedded into a half-plane that is its boundary does not intersect the $x$-axis.

In Figures 2 and 3, $\tau_1/\mu_0=-1$ and $\tau_1^\infty=-2$, the parameter $m$ is fixed, while the parameter
$\Gk$ has different values. In Fig.2, $\Gk\in(0,1)$ and $\Gk>1$ in Fig. 3. If $\Gk\to 1^\pm$, then the dimensionless size
of the inclusion is growing, and in order to prevent the intersection of the inclusion boundary with the $x$-axis,
it is necessary to increase the parameter $m$. 

In addition to the singular case $\Gk=1$, there is another 
singular case, $\tau_1=\tau_1^\infty$, when the solution to the problem does not exist.
In Figures 4 and 5, we show the inclusion profile when $\Gk=0.5$,
$\tau_1^\infty/\mu_0=-2$ for $m=2$ (Fig. 4), and $m=5$ (Fig. 5) for some values of the parameter $\tau_1/\mu_0$.
When the values of the stress $\tau_1$ inside the inclusion  
are  close to the values of the stress at infinity $\tau_1^\infty$, then the inclusion size decreases. 
 What is different between the cases $\tau_1/\tau_1^\infty<1$ and 
$\tau_1/\tau_1^\infty>1$ is that the function $\Go(\Gz)$ maps the parametric $\Gz$-plane onto the lower half-plane
in Fig. 4 and the upper half-plane in Fig. 5. Another difference is that when $\tau_1\to\tau_1^\infty$ and
 $|\tau_1/\tau_1^\infty|>1$ 
the inclusion approaches the half-plane boundary, while when  $|\tau_1/\tau_1^\infty|<1$ it drifts away from the boundary
of the external body, At the same time when  $|\tau_1/\tau_1^\infty|\to 0$, the inclusion is growing, becomes
closer to the $x$-axis, and we need to increase the parameter $m$ to keep the inclusion inside the external body.

\setcounter{equation}{0}

\section{Two scalar Riemann-Hilbert problems associated with the problem of $n$ inclusions}

In this section we aim to generalize the method by considering  $n\ge 2$ uniformly stressed inclusions $D_j$ 
embedded into a half-plane ${\Bbb R}_+^2=\{|x_1|<\infty,x_2\ge0\}$.
Suppose that the shear moduli of the inclusions $D_j$ and the surrounding infinite body $D_0={\Bbb R}_+^2\setminus\cup_{j=1}^n D_j$ are $\mu_j$ and $\mu_0$, respectively. At infinity, $\tau_{13}=\tau_1^\infty$, $\tau_{23}=0$.
The boundary of the half-plane is free of traction, $\tau_{23}=0$, and inside the inclusions, the stresses are
constant, $\tau_{13}=\tau^j_1$, $\tau_{23}=\tau^j_2$, $(x_1,x_2)\in D_j$, $j=1,2,\ldots, n$. As in Section 2, the interface conditions of ideal contact 
\beq
w_0=w_j, \quad \mu_0\fr{\Md w_0}{\Md \nu}=\mu_j\fr{\Md w_j}{\Md \nu}, \quad x\in L_j,\quad j=1,2,\ldots,n,
\label{6.1}
\eeq
hold on the boundary $L_j$ of an inclusion $D_j$.

\subsection{Case $\tau_1^j=\tau_1$, $\tau_2^j=\tau_2$, $j=1,2,\ldots,n$}

By introducing a function $f(z)$ by (\ref{2.5}) we rewrite the interface conditions (\ref{6.1}) as
\beq
f(z)=\fr{1}{\Gl_j}\R\fr{\bar\tau z}{\mu_0} +a_j+ib_j, \quad z\in L_j, \quad j=1,2,\ldots,n,
 \label{6.2}
 \eeq
 where $\Gl_j=\Gk_j/(1-\Gk_j)$, $\Gk_j=\mu_j/\mu_0$, and $a_j$ and $b_j$ are arbitrary real constants.

Let $z=\Go(\Gz)$ be a conformal map $\CD\to D_0$ from a parametric $\Gz$-plane cut along
segments $l_0=[m,\infty)$ and  $l_j=[k_{2j-1},k_{2j}]$, $j=1,2,\ldots,n$, where
$m>1$, $k_1=0$, $k_2=1$, $0>k_4>k_3>k_6>k_5>\ldots>k_{2n}>k_{2n-1}$.
The function $\Go(\Gz)$ maps the two-sided  finite segments $l_j$ into the contours $L_j$ and the two-sided semi-infinite
contour $l_0$ into the $x_1$-axis of the physical plane. When there are two inclusions in the half-plane and the domain 
is triply-connected, such a map always
exists. In the case  when $n\ge 3$, we assume that the inclusions are arranged such that all preimages $l_j$
of their boundaries $L_j$ lie in the real axis.  As in Section 3, we need to determine two functions,
$F(\Gz)=f(\Go(\Gz))$ and $\Go(\Gz)$ analytic in the $(n+1)$-connected slit domain $\CD$
based on the boundary conditions
$$
\I F(\Gz)=b_j, \quad \Gz\in l_1,\quad
\R \fr{\bar\tau\Go(\Gz)}{\mu_0}=\Gl_j[\R F(\Gz)-a_j], \quad \Gz\in l_j.\quad j=1,\ldots,n,
$$
\beq
\I \Go(\Gz)=0,\quad  \I F(\Gz)=b_0-\I\fr{\bar\tau\Go(\Gz)}{\mu_0}, \quad \Gz\in l_0.
\label{6.3}
\eeq
At infinity the functions $F(\Gz)$ and $\Go(\Gz)$ have to satisfy the conditions (\ref{3.5}).
 
 Let $\CR$ be a genus-$n$ hyperelliptic surface of the algebraic function 
 \beq
 u^2=p(\Gz), \quad p(\Gz)=(\Gz-m)\prod_{j=1}^{2n}(\Gz-k_j).
 \label{6.4}
 \eeq
Introduce two functions $\GF_1(\Gz,u)$ and $\GF_2(\Gz,u)$  on the 
surface $\CR$. The first function is defined as in (\ref{3.7}), while the second function is
\beq
\GF_2(\Gz,u)=\left\{\begin{array}{cc}
i\mu_0^{-1}\bar\tau\Go(\Gz), & (\Gz,u)\in\CD^+\subset\CR,\\
-i\mu_0^{-1}\tau\ov{\Go(\bar\Gz)}, & (\Gz,u)\in\CD^-\subset\CR.\\
\end{array}
\right.
\label{6.5}
\eeq
These functions satisfy the symmetry condition (\ref{3.7'}), the conditions (\ref{3.12}) at infinity,
analytic on the Riemann surface $\CR\setminus\CL$,
$\CL=\cup_{j=0}^n l_j$, and H\"older-continuous up to the two-sided contour $\CL$. On the contour $\CL$,
they satisfy the vector Riemann-Hilbert boundary condition (\ref{3.9}), where  $G(\Gx,v)$ is a piece-wise constant matrix
given by
\beq
G(\Gx,v)=\left(\begin{array}{cc}
1 & 0\\
2i\Gl_j & 1\\
\end{array}
\right), \quad (\Gx,v)\in l_j, \quad 
G(\Gx,v)=\left(\begin{array}{cc}
1 & i(1-\bar\tau/\tau)\\
0 & -\bar\tau/\tau\\
\end{array}\right), \quad (\Gx,v)\in l_0,
\label{6.6}
\eeq
and $g(\Gx,v)$ is a piece-wise constant vector
\beq
g(\Gx,v)=\left(\begin{array}{c}
2ib_j\\
-2i\Gl_j(a_j-ib_j) \\
\end{array}
\right), \quad (\Gx,v)\in l_j, \quad 
g(\Gx,v)=\left(\begin{array}{c}
2ib_0 \\
0 \\
\end{array}\right), \quad (\Gx,v)\in l_0.
\label{6.7}
\eeq 
When $\tau_2=0$, the vector Riemann-Hilbert problem is equivalent to two scalar problems
$$
\GF_1^+(\Gx,v)-\GF_1^-(\Gx,v)=2i b_j, \quad (\Gx,v)\in l_j, \quad j=0,1,\ldots,n,
$$
$$
\GF_2^+(\Gx,v)+\GF_2^-(\Gx,v)=0, \quad (\Gx,v)\in l_0,
$$
\beq
\GF_2^+(\Gx,v)-\GF_2^-(\Gx,v)=2i\Gl_j[\R\GF_1^+(\Gx,v)-a_j], \quad (\Gx,v)\in l_j,\quad  j=1,\ldots,n.
\label{6.8}
\eeq

\subsection{Case $\tau_1^j/\mu_j=\nu_1$, $\tau_2^j/\mu_j=\nu_2 $, $j=1,2,\ldots,n$}

In this case we denote $\nu=\nu_1+i\nu_2$ and rewrite the interface conditions
in the form
\beq
f(z)=i(\Gk_j-1)\I(\bar\nu z) +a_j+ib_j, \quad z\in L_j, \quad j=1,2,\ldots,n,
 \label{6.9}
 \eeq
The counterpart of the boundary conditions (\ref{6.3}) on the contours $l_j$ is
$$
\R F(\Gz)=a_j, \quad 
(\Gk_j-1)\I [\bar\nu \Go(\Gz)]=\I F(\Gz)-b_j, \quad \Gz\in l_j.\quad j=1,\ldots,n,
$$
\beq
\I \Go(\Gz)=0,\quad 
 \I F(\Gz)=b_0-
 \I [\bar\nu\Go(\Gz)], \quad \Gz\in l_0.
\label{6.10}
\eeq
To reduce this problem to a vector Riemann-Hilbert problem on the surface $\CR$ we
define the function $\GF_1(\Gz,u)$ and  $\GF_1(\Gz,u)$ in a way different from the case of Section 6.1. We put
$$
\GF_1(\Gz,u)=\left\{\begin{array}{cc}
iF(\Gz), & (\Gz,u)\in\CD^+\subset\CR,\\
-i\ov{F(\bar\Gz)}, & (\Gz,u)\in\CD^-\subset\CR,\\
\end{array}
\right.
$$
\beq
\GF_2(\Gz,u)=\left\{\begin{array}{cc}
\bar\nu\Go(\Gz), & (\Gz,u)\in\CD^+\subset\CR,\\
\nu\ov{\Go(\bar\Gz)}, & (\Gz,u)\in\CD^-\subset\CR.\\
\end{array}
\right.
\label{6.11}
\eeq
These functions solve the vector Riemann-Hilbert problem (\ref{3.9}) on the contour $\CL$
with the matrix coefficient $G(\Gx,v)$ and the right-hand side $g(\Gx,v)$ having the form
$$
G(\Gx,v)=\left(\begin{array}{cc}
1 & 0\\
-\fr{2i}{\Gk_j-1} & 1\\
\end{array}
\right), \; (\Gx,v)\in l_j, \quad 
G(\Gx,v)=\left(\begin{array}{cc}
-1 & i(1-\fr{\bar\nu}{\nu})\\
0 & \fr{\bar\nu}{\nu}\\
\end{array}\right), \; (\Gx,v)\in l_0,
$$
\beq
g(\Gx,v)=\left(\begin{array}{c}
2ia_j\\
-\fr{2i}{\Gk_j-1}(b_j+ia_j) \\
\end{array}
\right), \; (\Gx,v)\in l_j, \quad 
g(\Gx,v)=\left(\begin{array}{c}
-2b_0 \\
0 \\
\end{array}\right), \; (\Gx,v)\in l_0.
\label{6.12}
\eeq
Again, the vector problem is decoupled if $\tau_2^j=0$, that is $\nu=\bar\nu=\tau_1^j/\mu_j$, $j=1,\ldots,n$.
We deduce the following two scalar Riemann-Hilbert problems on the genus-$n$ hyperelliptic surface $\CR$:
$$
\GF_1^+(\Gx,v)-\GF_1^-(\Gx,v)=2i a_j, \quad (\Gx,v)\in l_j, \quad j=1,\ldots,n
$$
$$
\GF_1^+(\Gx,v)+\GF_1^-(\Gx,v)=-2b_0, \quad (\Gx,v)\in l_0,
$$
$$
\GF_2^+(\Gx,v)-\GF_2^-(\Gx,v)=-\fr{2i}{\Gk_j-1}[\R\GF_1^+(\Gx,v)+b_j], \quad (\Gx,v)\in l_j,\quad  j=1,\ldots,n,
$$
\beq
\GF_2^+(\Gx,v)-\GF_2^-(\Gx,v)=0,  \quad (\Gx,v)\in l_0.
\label{6.13}
\eeq
In both Riemann-Hilbert problems,  (\ref{6.8})  and   (\ref{6.13}),  the functions $\GF_1(\Gz,u)$ and  $\GF_2(\Gz,u)$ satisfy the symmetry condition (\ref{3.7'}) and the conditions at 
infinity (\ref{3.12}), and it is required to solve the problem of factorization. 

For the solution of both Riemann-Hilbert problems one need a genus-$n$ hyperelliptic analogue of the elliptic kernel (\ref{4.4}). 
It has the form
 \beq
dV((\Gx,v),(\Gz,u))=\fr12\left[1+\fr{u}{v}\prod_{j=0}^{n}\fr{\Gx-\Gx_j}{\Gz-\Gx_j}\right]
\left(\fr{1}{\Gx-\Gz}-\fr{1}{\Gx-\Gx_0}\right)d\Gx.
\label{6.14}
\eeq 
Here, $\Gx_0,\ldots,\Gx_n$ are arbitrary fixed not necessarily distinct points in the real axis not falling on the 
 contour $\CL$ of the Riemann-Hilbert problem. 
The solution of the Riemann-Hilbert problems (\ref{6.8})  and (\ref{6.13}) is beyond the scope of this paper and not presented here.

\section{Conclusion}

In the previous sections a closed-form solution has been derived to the inverse problem of antiplane shear of 
an elastic finite domain $D_1$ embedded into an elastic  half-plane ${\Bbb R}_+^2$. In this model, 
the boundary of the half-plane $\{||x_1|<\infty, x_2=0\}$
is kept   free of traction $\tau_{23}$, the conditions of ideal contact on the interface hold, and  the stress 
field inside the inclusion is uniform, while the shape of the inclusion is to be recovered.
This harmonic problem is well suited for the method of conformal mappings.
It has been shown that a map from a parametric plane cut along  the segments $[0,1]$   and $[m,\infty)$
 onto the physical doubly-connected domain  ${\Bbb R}_+^2\setminus D_1$ can be recovered by solving a 
vector Riemann-Hilbert problem with a piece-wise constant matrix coefficient on an elliptic
surface $\CR$. 
Under the assumption
that $\tau_{23}=0$ inside the inclusion $D_1$, the vector 
problem has been decoupled into two scalar Riemann-Hilbert problems on two slits on the surface $\CR$.
The solution has been determined by proposing and using an analogue of the Cauchy kernel on an elliptic surface and
solving a Jacobi inversion problem associated with the factorization problem, a part of the solution procedure. 
Not counting  the scaling parameter, the solutions to the Riemann-Hilbert problems 
comprise a four-parametric family of conformal mappings which possess three model parameters
$\Gk=\mu_1/\mu_0$, $\tau_1/\mu_0$, and $\tau_1^\infty/\mu_0$ and one geometric parameter $m$.
Here, $\mu_1$ and $\mu_0$ are the shear moduli of the inclusion and the half-plane, $\tau_1$ is the constant stress
$\tau_{13}$ inside the inclusion, and $\tau_1^\infty$ is the stress $\tau_{13}$ applied at infinity. 
Numerical implementation of this method has shown  that
the inclusion is not symmetric with respect to any line parallel to the $x_1$-axis. At the same time,
its shape resembles an ellipse. 
 
The method can be generalized to the problem  of $n$ inclusions $D_j$ ($j=1,\ldots,n)$  in a half-plane
when uniform stresses inside all the inclusions satisfy one of the conditions,  
(i)  $\tau_{13}^j=\tau_1$, $\tau_{23}^j=\tau_2$, or
(ii)  $\tau_{13}^j/\mu_j=\nu_1$, $\tau_{23}^j/\mu_j=\nu_2$.
The conformal map, if available, recoveres  the whole family of inclusions in the case $n=2$
and a part of the family if $n\ge 3$ that is the set of those inclusions whose preimages lie in the real axis
of the parametric plane. If $\tau_2\ne 0$ ($\nu_2\ne 0$),  one needs to deal with a vector Riemann-Hilbert problem on a genus-$n$
hyperelliptic surface. In the particular case $\tau_2=0$ ($\nu_2=0$), the problem is decoupled, and the corresponding
scalar Riemann-Hilbert problems admit a solution by singular integrals. One of the most difficult issues 
in the solution procedure is a genus-$n$ Jacobi inversion problem whose solution is determined through
the $n$ zeros of the associated Riemann Theta function \cite{zve}, \cite{as}.

\vspace{.1in}

{\bf Data accessibility.}
 No software generated data were created during this study.
 
{\bf Competing interests.} I declare I have no competing interests.

{\bf Funding statement.} This research received no specific grant from any funding agency in the public, commercial or 
not-for-profit sectors.


\begin{thebibliography}{99}

\bibitem{gus}
Gustafsson B, Vasilev A. 2006 \textit{Conformal and potential analysis in Hele-Shaw cells}. Advances in Mathematical Fluid Mechanics. Basel: Birkhäuser Verlag.

\bibitem{cro}
Crowdy DG. 2006 Exact solutions to the unsteady two-phase Hele-Shaw problem. 
\textit{Quart. J. Mech. Appl. Math.}   \textbf{59}, 475-485. 

\bibitem{antsil}
Antipov YA, Silvestrov VV. 2007
Method of Riemann surfaces in the study of supercavitating flow around two hydrofoils in a channel. 
\textit{Physica D}   \textbf{235},  72-81. 

\bibitem{chr} Christopher TW,  Llewellyn Smith SG. 2021
Hollow vortex in a corner, \textit{J. Fluid Mech.}  \textbf{908},  R2 1-12.

\bibitem{antzem} Antipov YA, Zemlyanova AY. 2021 
Sadovskii vortex in a wedge and the associated Riemann-Hilbert problem on a torus.  arXiv:2010.08118v2, 23 pages.

\bibitem{esh}
Eshelby JD. 1957 The determination of the elastic field of an ellipsoidal inclusion, and related problems.
 \textit{Proc. Roy. Soc. London A.}  \textbf{ 241},  376-396.

\bibitem{che}
Cherepanov GP. 1974  Inverse problems of the plane theory of elasticity.  \textit{J. Appl. Math. Mech. PMM}  \textbf{38}, 915-931. 

\bibitem{vig}
 Vigdergauz SB.  1976 Integral equation of the inverse problem of the plane theory of elasticity.
\textit{J. Appl. Math. Mech. PMM} \textbf{40}, 518-522.

\bibitem{ant1}
Antipov YA. 2018 Slit maps in the study of equal-strength cavities in $n$-connected elastic planar domains. 
\textit{ SIAM J. Appl. Math.}  \textbf{78}, 320–342.

\bibitem{mar}
Marshall JS. 2019 On sets of multiple equally strong holes in an infinite elastic plate: parameterization and existence. \textit{SIAM J. Appl. Math.} \textbf{79}, 2288–2312. 

\bibitem{obn}
Obnosov Y, Zulkarnyaev A.  2019 Nonlinear mixed Cherepanov boundary value problem.  
\textit{Complex Var. Elliptic Equ.}    \textbf{64}, 979–996.

\bibitem{che2}
Cherepanov GP. 1962
Inverse elastic-plastic problem for antiplane strain. \textit{J. Appl. Math. Mech. PMM}
\textbf{26}, 1743-1748.

\bibitem{kan}
Kang H, Kim E,  Milton GW. 2008  Inclusion pairs satisfying Eshelby's uniformity property.
\textit{SIAM, J. Appl. Math}  \textbf{69},  577-595.

\bibitem{wan1} Wang X. 2012 Uniform fields inside two
non-elliptical inclusions.  \textit{Math. Mech. Solids} \textbf{17}, 736-761.

\bibitem{wan2} Wang X, Yang P, Schiavone P. 2020
Uniform stresses inside a non-elliptical inhomogeneity and a nearby half-plane with
locally wavy interface. \textit{Z. Angew. Math. Phys.}  \textbf{ 71}, paper 58, 1-11.

\bibitem{liu}
Liu LP.  2008 Solutions to the Eshelby conjectures.   \textit{Proc. Roy. Soc. London A.}  \textbf{ 464}, 573-594.

\bibitem{dai}
Dai M, Ru CQ,  Gao CF. 2017  Uniform strain fields inside multiple inclusions in an elastic infinite plane under
anti-plane shear.  \textit{Math. Mech. Solids.} \textbf{ 17}, 114-128.

\bibitem{ant2}
Antipov YA. 2020 Method of Riemann surfaces for an inverse antiplane problem in an n-connected domain.
\textit{ Complex Var. Elliptic Equ.}   \textbf{65},  455–480

\bibitem{ant3}
 Antipov YA. 2019 Method of automorphic functions for an inverse problem of antiplane elasticity. 
 \textit{Quart. J. Mech. Appl. Math.}  \textbf{72}, 213–234.

\bibitem{az}  Antipov YA, Zemlyanova AY. 2009 Motion of a yawed supercavitating wedge beneath a free surface. 
 \textit{SIAM J. Appl. Math.}  \textbf{70},  923-948.

\bibitem{zve} Zverovich EI. 1971 Boundary value problems in the theory of analytic functions in H\"older classes on Riemann surfaces.  \textit{Russian Math. Surveys}  \textbf{26}, 117-192.


\bibitem{as} 
Antipov YA, Silvestrov VV. 2006  
Electromagnetic scattering from an anisotropic impedance half plane
at oblique incidence: the exact solution. \textit{Quart. J. Mech. Appl. Math.}  \textbf{259},  211-251.



\end{thebibliography}
\end{document}